\documentclass{article}  

\usepackage{amssymb}
\usepackage{amsmath}
\usepackage{amsfonts}
\usepackage{amscd}
\usepackage{amsthm}
\usepackage{hyperref}

\newtheorem{definition}{Definition}

\newtheorem{question}{Question}

\newcommand{\Q}{\mathbb{Q}}

\title{Low Degree Places on the Modular Curve $X_1(N)$}

\author{Mark van Hoeij\thanks{Supported by NSF grants 1017880 and 1319547.}  \\
Florida State University, Tallahassee, FL 32306-3027, USA \\
hoeij@math.fsu.edu}

\begin{document}
\maketitle
\noindent Let $\alpha := 2 \cos(2 \pi/7)$ and let $\tau$ be the golden
ratio, they are solutions of
\[ \alpha^3 + \alpha^2 - 2\alpha - 1 =0, \ \ \tau^2-\tau-1 = 0. \]
Let
\[ b := (6\tau-3)\alpha^2 + (14\tau - 8)\alpha + 5\tau - 3 \]
and
\[ c := \tau  \alpha^2 + 2\tau \alpha + 1. \]
Let $E_{b,c}$ be the elliptic curve
\begin{equation}
	E_{b,c} := y^2 + (1-c)xy - b  y = x^3 - b  x^2   \label{Ebc}
\end{equation}
and
\[ P := (x=0, y=0). \]
Then $P$ is a point on $E_{b,c}$ of order 37.  The pair $(E_{b,c}, \, P)$
corresponds to a point on the modular curve $X_1(37)$.
The interest in this example lies in the fact that this point is defined
over a number field $\Q(\alpha,\tau)$ of degree 6, whereas the $\Q$-gonality
of $X_1(37)$ is 18.  This short note is motivated by the following
question: Are there only finitely many points on $X_1(N)$ with degree less than the $\Q$-gonality,
and if so, can they all be found with a finite computation? \\[10pt]

% The purpose of this short note is this and related questions, but first some notation.
\noindent {\bf Notations:} Let $C_N$ denote the function field of $X_1(N)$ over $\Q$.
We can write $C_N = \Q(x)[y]/(F_N)$ where $F_N$ is an explicit equation given by Sutherland
at: \url{http://math.mit.edu/~drew/X1_altcurves.html}. \ \ 
	% For $N=10,11,\ldots$, let $F_N \in \Q[x,y]$ denote
	% the explicit equation given by Sutherland
	% for the modular curve $X_1(N)$, in the files: \\
% \verb+http://math.mit.edu/~drew/FFFc10.txt+, \ \ \verb+FFFc11.txt+, \ldots. \\
% Denote $C_N := \Q(x)[y]/(F_N)$, the function field of $X_1(N)$ over $\Q$.
The $\Q$-gonality of $X_1(N)$ is
\[ {\rm gon}(N) := {\rm min}\{ {\rm degree}(f) \, | \, f \in C_N - \Q\}. \]
In a joint work with Maarten Derickx, ${\rm gon}(N)$ has been computed for all $N \leq 40$, see  arXiv:1307.5719.

Denote ${\rm Places}(N)$ as the set of discrete valuations on
$C_N$ over $\Q$.
	%  i.e. the set of all onto maps
	% \[ v: C_N \rightarrow \Z \bigcup \{\infty\} \]
	% with (i) $v(f_1 f_2) = v(f_1) + v(f_2)$,  \ \
	% (ii) $v(f_1 + f_2) \geq {\rm min}(v(f_1),v(f_2))$
	% (iii) $v(f)=\infty  \Longleftrightarrow f=0$,
	% and $v(q)=0$ for all $q \in \Q - \{0\}$.
For each place $v$, denote $\Q(v)$ as the residue-field of $v$, and
${\rm deg}(v) := [\Q(v): \Q]$.
Consider the functions $r,s,b,c \in C_N$ written on Sutherland's webpage as:
\[ r := \frac{x^2y-xy+y-1}{x^2y-x}, \ \  s := \frac{xy-y+1}{xy}, \ \  b  := rs(r-1), \ \ c := s(r-1).  \]
Then let $j \in C_N$ be the $j$-invariant of the curve $E_{b,c}$ from equation~(\ref{Ebc}).
	% Then $j/1728$ is a Belyi map, the natural map from $X_1(N)$ to $X_1(1)$.
The {\em cusps} are the places where $E_{b,c}$ degenerates,
% i.e., the places where $j$ has a pole,
\[ {\rm Cusps}(N) := \{v \in {\rm Places}(N) \, | \, v(j) < 0 \}. \]
If $v$ is not a cusp, then $v$ corresponds to an elliptic curve $E_{b,c}$ over $\Q(v)$ with a point $P$ of order $N$. \\[5pt]

% In a joint work with Maarten Derickx, ${\rm gon}(N)$ has been computed for all $N \leq 40$, see  arXiv:1307.5719.
% and upper bounds were given for 
% The paper is in progress, but the gonality data is available: \\
% \verb+http://www.math.fsu.edu/~hoeij/files/X1N/gonality+ \\
	% For each $N \leq 40$ we computed an upper bound by computing an explicit $f \in C_N$, and a lower bound by
	% computing the $\Z/(p)$-gonality for the smallest prime $p$ that does not divide $N$.  The lower and upper
	% bounds match, and hence the above URL gives the exact $\Q$-gonality for $N \leq 40$.
% This URL also enumerates all decompositions of $X_1(N) \rightarrow X_1(1)$ (i.e. decompositions
% of the function $j$).  It also lists upper bounds for ${\rm gon}(N)$ for $N \leq 100$
% (these upper bounds are likely sharp, but that is only proven for $N \leq 40$). \\[5pt]

\begin{definition}
Define the set of {\em \bf low-degree places} as
\[
        {\rm LDP}(N) := \{ v \in {\rm Places}(N) - {\rm Cusps}(N) \ | \ {\rm deg}(v) < {\rm gon}(N) \}.
\]
Let $f \in C_N - \Q$.  The {\rm support} of $f$ is ${\rm Supp}(f) := \{ v \in {\rm Places}(N) \, | \, v(f) \neq 0\}$.
A function $f$ is called a {\em modular unit} if ${\rm Supp}(f) \subseteq {\rm Cusps}(N)$. \\
A {\em cusp-neighbor} is a place $v \in {\rm Places}(N) - {\rm Cusps}(N)$ for which there exists $f \in C_N$
with $v(f) = 1$ and ${\rm Supp}(f) \subseteq {\rm Cusps}(N) \bigcup \{ v \}$.
We are only interested in low-degree cusp-neighbors, i.e. ${\rm deg}(v) < {\rm gon}(N)$.
\end{definition}

The example on page 1 shows ${\rm LDP}(37) \neq \emptyset$.  The main question is:
Is ${\rm LDP}(N)$ always finite, and can it be computed in a finite number of steps?

For any given $N$, the set of low-degree cusp-neighbors is finite, and can be computed in a finite number of steps.
The question, however, is if this would produce all low-degree places:
\begin{question} \label{QuestionNeighbor} Is every low-degree place a cusp-neighbor?
\end{question}

The author has computed a number of low-degree places by
computing Riemann-Roch spaces $L(D)$ for suitably chosen divisors $D$ with ${\rm Supp}(D) \subseteq {\rm Cusps}(N)$.
	% In order to decide which $D$'s to use in characteristic $0$, first compute the Picard group in characteristic $p$.
Combining this with $L(D)$-computations over a finite field\footnote{To decide which $D$'s to use,
consider all divisors $D'$ over $\mathbb{F}_p$
with $D' \geq 0$, ${\rm deg}(D') < {\rm gon}(N)$, and
${\rm Supp}(D') \bigcap {\rm Cusps}(N) = \emptyset$.
For each such $D'$, compute, if it exists, a $D$ with ${\rm Supp}(D) \subseteq {\rm Cusps}(N)$
for which $D-D'$ is principal. Lift such $D$ to characteristic 0.}
the author estimates that it should be possible to find all low-degree cusp-neighbors on a computer for $N \leq 40$.
This strategy has not yet been systematically implemented, but the author has tested it to see if it can produce
explicit examples, which it does:
% When more examples are computed they will be posted on: \verb+http://www.math.fsu.edu/~hoeij/files/X1N/+

\newpage
\noindent {\bf \large Examples:} \\[5pt]
The $\Q$-gonality of $X_1(29)$ is 11.
A low-degree place of $X_1(29)$ is given below, in the form $(E_{b,c}, \, P)$,
for three non-isomorphic number fields.
Additional places over the same fields can be found by taking multiples of $P$.
	% multiplying $P$ with an element of $\{1,\ldots,14\}$.
To decide whether or not $X_1(29)$ has more low-degree places, one needs an answer to Question 1, and, a
systematic implementation. 
	%  For each, one element of $X_1(29)$
	% is given; additional elements can be found by computing a multiple of the point $P$.
	% To decide whether or not there are more, one would need to answer Question 1, and write a
	% systematic implementation.
The examples can be copied from the file LowDegreePlaces at
\url{http://www.math.fsu.edu/~hoeij/files/X1N}.  \\[5pt]
Let $a$ be a solution of $a^9-a^8-5a^7+5a^6+7a^5-8a^4-2a^3+4a^2-a-1 = 0.$ \\
$b := a^7-5a^6-3a^5+22a^4-3a^3-28a^2+9a+8$,  \\
$c := -2a^8+2a^7+10a^6-10a^5-16a^4+15a^3+10a^2-6a-3$. \\
Then $E_{b,c}$ has a point $P = (0,0)$ of order $N=29$.  \\[6pt]
	% Multiplying $P$ by elements of $\{1,\ldots,14\}$ produces
	% more elements of ${\rm Places}(29)$ with residue-field $\Q(a)$.
	%
	% \label{examples}
	% $a^10-2*a^7+a^5+2*a^4+3*a^3-3*a^2-2*a+1 = 0$, 
	% $b = (-145*a^9-14*a^8-66*a^7+114*a^6+175*a^5-137*a^4-211*a^3-324*a^2-85*a+92)/43$, \\
	% $c = a^9+a^7-a^6-a^5+a^4+a^3+3*a^2+a-1$.  Then $E_{b,c}$ has a point of order $N=29$. Multiplying $P$
	% by an element of $\{1,\ldots,14\}$ produces 14 pairs $(b,c)$, i.e. 14 elements of ${\rm Places}(29)$ with residue-field $\Q(a)$.
\noindent Let $a^{10}-2a^7+a^5+2a^4+3a^3-3a^2-2a+1 = 0$, \\
$b := (-145a^9-14a^8-66a^7+114a^6+175a^5-137a^4-211a^3-324a^2-85a+92)/43$, \\
$c := a^9+a^7-a^6-a^5+a^4+a^3+3a^2+a-1$. \\[6pt]
	% Then $E_{b,c}$ has a point $P$ of order $N=29$. Multiplying $P$
	% by an element of $\{1,\ldots,14\}$ produces 14 elements of ${\rm Places}(29)$ with residue-field $\Q(a)$. \\
	%
	% \noindent Same degree but different number field: \\
	% $a^10-2*a^9+2*a^8-5*a^7+7*a^6-4*a^5+4*a^4-8*a^3+5*a^2-1 = 0$, \\
	% $b = (-23*a^9+21*a^8+19*a^7-71*a^6+95*a^5-159*a^4+178*a^3-78*a^2-10*a+6)/97$, \\
	% $c = (60*a^9-59*a^8+39*a^7-207*a^6+98*a^5+31*a^4+8*a^3-58*a^2+5*a-3)/97$,  again $X_1(29)$.  \\[5pt]
	% \noindent Same degree but a different number field, again for $X_1(29)$: \\
\noindent Let $a^{10}-2a^9+2a^8-5a^7+7a^6-4a^5+4a^4-8a^3+5a^2-1 = 0$, \\
$b := (-23a^9+21a^8+19a^7-71a^6+95a^5-159a^4+178a^3-78a^2-10a+6)/97$, \\
$c := (60a^9-59a^8+39a^7-207a^6+98a^5+31a^4+8a^3-58a^2+5a-3)/97$.  \\[10pt]
	% Combined with the example on page 1, this gives 3 number fields over which a low-degree place can
	% be found.  To decide whether or not there are more, one would need to answer Question 1, and write a
	% systematic implementation. \\
	%
	% \noindent In the next sequence of examples, we take $N=31$, which has $\Q$-gonality 12. We give examples
	% for 5 non-isomorphic fields, one of which has degree 9, which well below the $\Q$-gonality. \\
	% $a^9-3*a^8+4*a^7-a^6-7*a^5+11*a^4-9*a^3+3*a^2-a+1 = 0$, \\
	% $b = (-90*a^8+228*a^7-335*a^6+205*a^5+183*a^4-246*a^3+229*a^2+49*a+66)/37$,
	% $c = (-48*a^8+129*a^7-117*a^6-51*a^5+364*a^4-368*a^3+169*a^2-a+50)/37$. \\[10pt]
\noindent The next examples are for $X_1(31)$. Five non-isomorphic fields
are given, one of which has degree 9. The $\Q$-gonality of $X_1(31)$ is 12. \\[5pt]
Let $a^9-3a^8+4a^7-a^6-7a^5+11a^4-9a^3+3a^2-a+1 = 0$, \\
$b := (-90a^8+228a^7-335a^6+205a^5+183a^4-246a^3+229a^2+49a+66)/37$, \\
$c := (-48a^8+129a^7-117a^6-51a^5+364a^4-368a^3+169a^2-a+50)/37$. \\[6pt]
	%
	%
	% \noindent $a^10+2*a^8-3*a^7+3*a^6-7*a^5+8*a^4-7*a^3+7*a^2-4*a+1 = 0$,
	% $b = 62*a^9-48*a^8+47*a^7-192*a^6+262*a^5-321*a^4+421*a^3-330*a^2+131*a-21$,
	% $c = 8*a^9+6*a^8+4*a^7-16*a^6+2*a^5-8*a^4+10*a^3+14*a^2-16*a+5$.  \\[5pt]
	%
	% \noindent $a^11-2*a^10-3*a^9+9*a^8-a^7-13*a^6+9*a^5+7*a^4-5*a^3-1 = 0$,
	% $b = -9*a^10+38*a^9-35*a^8-59*a^7+131*a^6-45*a^5-81*a^4+59*a^3+6*a^2+2*a+14$,
	% $c = 3*a^9-7*a^8-3*a^7+18*a^6-9*a^5-13*a^4+10*a^3+9*a^2+2*a$. \\[5pt]
	%
	% \noindent $a^11-4*a^10+9*a^9-15*a^8+21*a^7-21*a^6+17*a^5-8*a^4+3*a^2-3*a+1 = 0$,
	% $b = -a^10-3*a^9-4*a^8-6*a^7-6*a^6-a^5-6*a^4+6*a^3-2*a^2+3*a+1$,
	% $c = -a^10+2*a^9-3*a^8+3*a^7-5*a^6-a^5-a^4-3*a^3+2*a^2-a+2$. \\[5pt]
	%
	% \noindent $a^11-a^10-4*a^9+7*a^8+4*a^7-9*a^6-5*a^5+2*a^4+8*a^3+2*a^2-3*a-1$,
	% $b = (-245*a^10+1414*a^9+3377*a^8-3908*a^7-7202*a^6+562*a^5+5683*a^4+5190*a^3-449*a^2-2406*a-678)/349$,
	% $c = (8*a^10-106*a^9-304*a^8+290*a^7+842*a^6-440*a^5-932*a^4+265*a^3+395*a^2-24*a-79)/349$.  \\[5pt]
	%
	%
\noindent Let $a^{10}+2a^8-3a^7+3a^6-7a^5+8a^4-7a^3+7a^2-4a+1 = 0$, \\
$b := 62a^9-48a^8+47a^7-192a^6+262a^5-321a^4+421a^3-330a^2+131a-21$, \\
$c := 8a^9+6a^8+4a^7-16a^6+2a^5-8a^4+10a^3+14a^2-16a+5$.  \\[6pt]
\noindent Let $a^{11}-2a^{10}-3a^9+9a^8-a^7-13a^6+9a^5+7a^4-5a^3-1 = 0$, \\
$b := -9a^{10}+38a^9-35a^8-59a^7+131a^6-45a^5-81a^4+59a^3+6a^2+2a+14$, \\
$c := 3a^9-7a^8-3a^7+18a^6-9a^5-13a^4+10a^3+9a^2+2a$. \\[6pt]
\noindent Let $a^{11}-4a^{10}+9a^9-15a^8+21a^7-21a^6+17a^5-8a^4+3a^2-3a+1 = 0$, \\
$b := -a^{10}-3a^9-4a^8-6a^7-6a^6-a^5-6a^4+6a^3-2a^2+3a+1$, \\
$c := -a^{10}+2a^9-3a^8+3a^7-5a^6-a^5-a^4-3a^3+2a^2-a+2$. \\[6pt]
\noindent Let $a^{11}-a^{10}-4a^9+7a^8+4a^7-9a^6-5a^5+2a^4+8a^3+2a^2-3a-1 = 0$, \\
$b := (-245a^{10}+1414a^9+3377a^8-3908a^7-7202a^6+562a^5+5683a^4+5190a^3-449a^2-2406a-678)/349$, \ \ \ \
$c := (8a^{10}-106a^9-304a^8+290a^7+842a^6-440a^5-932a^4+265a^3+395a^2-24a-79)/349$.  \\[10pt]
	%
	%
% \noindent We end with one more question (note: this holds for $N \leq 40$):
% \begin{question} Must $C_N$ have a modular unit of degree ${\rm gon}(N)$?
% \end{question}
	% There need not exist a function $f$ of degree ${\rm gon}(N)$ whose support
	% consists only of rational cusps (the first example of that is $N=36$). \\[5pt]
	%
\noindent The author would like to thank Maarten Derickx and Andrew Sutherland for fruitful discussions.
\newpage

\section{Degree table}
\newcommand{\ems}{\em \scriptsize}
\newcommand{\ns}{\mbox{}\hspace{-3pt}\mbox{}}
Aug 2013 update: The preceeding part of this preprint was written in February 2012. A small computation (just a few $N$'s) was done
to test if the approach was effective.
Here we extend the search to $N \leq 60$.
	%   (to get a better idea of the growth of the degrees, we plan to extend the
	%   computation to significantly larger $N$'s, but this requires implementing some algorithmic improvements).
We start with some notation regarding the table below.

A degree written in this font, 5,  indicates that we found an explicit function in $C_N$ of that degree.
The notation 4$^+$ (e.g. for $N=17$) means we have functions of degrees 4,5,6,7
implying that $C_N$ has functions of any degree $\geq 4$.

Degrees written in this font, {\ems 6},  indicate we found
non-cuspidal places $v$ on $X_1(N)$ of that degree. For example, the table indicates
that we found no functions of degrees 6,7 for $N=25$ % ($C_{25}$ has no functions of those degrees)
but we did find places of those degrees.

The notation {\ems 7}$^+$ indicates that $X_1(N)$ has non-cuspidal places of any degree $\geq 7$.  For example, the table
indicates that our file LowDegreePlaces on \url{http://www.math.fsu.edu/~hoeij/files/X1N}
for $N=33$ contains non-cuspidal places of degrees {\ems 7,8,9,11}  (for $N=33$ we do not store places
of degrees 10 and 12$^+$ because the explicit functions we have for those degrees can be used to construct arbitrarily many
places of those degrees).

% Another example is $N=37$, the table indicates the existence of places of degrees {\ems 6,10,12}$^+$
% and functions of degrees 18$^+$.
% It can be proven that $Q(X_1(37))$ has no functions of degree $<18$, but this is a significant amount of work, see arXiv:1307.5719.

At the moment, we have no proof that the table is complete.
However, with recent ideas from Maarten Derickx, it may be possible to perform a provably complete search for $N \leq 40$, $N \neq 37$.
\vspace{20pt}

\begin{tabular}{|c|l|c|l|c|l|}
\hline
$N$ 	& degrees 		& $N$	& degrees 				& $N$	& degrees					\\ \hline
\ns 1--10\ns& 1$^+$ 		& 29	& {\ems 9,10}, 11$^+$			& 45	& {\ems 10,12,14}$^+$, 18, 20$^+$		\\ \hline
11	& 2$^+$ 		& 30	& {\ems 5}, 6$^+$			& 46	& {\ems 14}$^+$, 19$^+$				\\ \hline
12	& 1$^+$ 		& 31	& {\ems 9}$^+$, 12$^+$			& 47	& {\ems 20}$^+$, 29$^+$				\\ \hline
\ns 13--16\ns& 2$^+$		& 32	& 8, {\ems 9}, 10$^+$			& 48	& {\ems 11,12,14}$^+$, 16, 18$^+$		\\ \hline
17	& 4$^+$			& 33	& {\ems 7}$^+$, 10, 12$^+$		& 49	& {\ems 14,19}, 21, {\ems 22}$^+$, 30$^+$	\\ \hline
18	& 2$^+$ 		& 34	& {\ems 8,9}, 10$^+$			& 50	& {\ems 10,12}, 15, {\ems 16}$^+$, 20, 22$^+$	\\ \hline
19	& 5$^+$			& 35	& {\ems 8,10}$^+$, 12, 14$^+$		& 51	& {\ems 15,18}$^+$, 24, 29$^+$			\\ \hline
20	& 3$^+$			& 36	& {\ems 7}, 8$^+$			& 52	& {\ems 16}$^+$, 21, 24$^+$			\\ \hline
21      & {\ems 3},\,4$^+$	& 37	& {\ems 6,10,12}$^+$, 18$^+$		& 53	& {\ems 22,25}$^+$, 37$^+$			\\ \hline
22	& 4$^+$			& 38	& {\ems 10}, 12$^+$			& 54	& {\ems 13,15}$^+$, 18, 20$^+$			\\ \hline
23	& 7$^+$			& 39	& {\ems 8--10,12}$^+$, 14, 16$^+$\ns	& 55	& {\ems 18,23}$^+$, 30, 34$^+$			\\ \hline
24	& 4$^+$ 		& 40	& {\ems 8}$^+$, 12, 14$^+$		& 56	& {\ems 18}$^+$, 24, 26, 28$^+$ 		\\ \hline
25	& 5,{\ems 6,7},\,8$^+$\ns& 41	& {\ems 14,17}$^+$, 22$^+$		& 57	& {\ems 12,16,18,19,21,22,24}$^+$, 30,36$^+\ns$\\ \hline
26	& 6$^+$			& 42	& {\ems 8}$^+$, 12$^+$			& 58	& {\ems 12,14,16,20}$^+$, 31$^+$		\\ \hline
27	& 6$^+$			& 43	& {\ems 12,14,15,17}$^+$,\,24$^+$\ns	& 59	& {\ems 31}$^+$, 46$^+$				\\ \hline
28	& {\ems 5},\,6$^+$	& 44	& {\ems 11}$^+$, 15$^+$			& 60	& {\ems 13,15}$^+$, 24, 26$^+$			\\ \hline
\end{tabular}
\newpage

June 2014 update: Added $N \in \{61,\ldots,80\}$.  For each $N$, the table indicates which function-degrees and which place-degrees
our search has produced (the search is probabilistic, there is no guarantee that this table is complete).
%
% if our search yielded a function or a place of degree $d$
%  The answer to question~\ref{QuestionNeighbor} is no. However, this does not imply that the degree-table is incomplete.
% Due to the probabilistic nature of the search, there is a non-trivial chance that some degree(s) is/are 
% missing, a chance that increases as a function of $N$ and the rank of the Jacobian.
	% So if any degree(s) are missing, it is more likely that this happens in this table than in the $N \leq 60$ table on the previous page.
\vspace{20pt}

\begin{tabular}{|c|l|c|l|}
\hline
$N$ & degrees  & $N$   & degrees  \\ \hline
61  &  {\ems 20,24,26,27,30,31,33}$^+$, 49$^+$\ns	&	71  &  {\ems 44,45,47}$^+$, 66$^+$			\\ \hline
62  &  {\ems 22}$^+$, 36$^+$				&	72  &  {\ems 22,24}$^+$, 32, 36, 40$^+$			\\ \hline
63  &  {\ems 18,20}$^+$, 36, 39, 41$^+$			&	73  &  {\ems 24,30,36,42,46,48}$^+$, 70$^+$		\\ \hline
64  &  {\ems 24}$^+$, 32, 36, 38$^+$			&	74  &  {\ems 18,20,29--31,34}$^+$, 51$^+$		\\ \hline
65  &  {\ems 20,24,26,28,30}$^+$, 42, 48$^+$		&	75  &  {\ems 25,31--33,35--37,39}$^+$, 40,45,50,55,60$^+$\ns\\ \hline
66  &  {\ems 16,19}$^+$, 30, 32, 35$^+$			&	76  &  {\ems 30,35}$^+$, 45, 48, 50, 52, 53, 54, 56$^+$	\\ \hline
67  &  {\ems 22,30,33,37,39,43}$^+$, 58$^+$		&	77  &  {\ems 40,48}$^+$, 60, 68, 72$^+$			\\ \hline
68  &  {\ems 26}$^+$, 36, 40, 42$^+$			&	78  &  {\ems 24,25,27,28,30}$^+$, 42, 48, 49, 51$^+$	\\ \hline
69  &  {\ems 28,29,32,34,36}$^+$, 44, 54$^+$		&	79  &  {\ems 26,42,51,54,57--59,61}$^+$, 82$^+$		\\ \hline
70  &  {\ems 20,24,26}$^+$, 36, 40, 42$^+$		&	80  &  {\ems 20,24,28,32,35}$^+$, 48, 54, 56$^+$	\\ \hline
\end{tabular}

\vspace{20pt}
The answer to question~\ref{QuestionNeighbor} turned out to be no.
%  for $N=37$ (the first $N$ for which the Jacobian of $X_1(N)$ over $\Q$ has positive rank).
\begin{question}
\label{Q2}
If $P$ is a low-degree place, must there be a modular unit $f$ over $\Q$ for which $f(P)  \in \Q$?
\end{question}
The motivation for this question is as follows:  For $N \leq 40$ we computed low-degree places via 
Riemann-Roch computations.  For $N>40$ these computations became slow, and so we switched to another
method:  Take random modular units $f$, and compute the roots of $f -c = 0$  (where $c$ is almost always $\pm 1$). \\
With this method we can reach higher values for $N$ than with the Riemann-Roch method. However,
it looks quite ad hoc, so we decided to test its effectiveness experimentally by applying
it to $N \leq 40$ as well.  It
turned out that every place $P$ we found for $N \leq 40$ with the Riemann-Roch method was again found
with the $f - c$ method.   So the $f - c$ method is surprisingly effective (at least, for values
of $N$ where we can compare with other methods).  The question is if there is a mathematical result
that could explain this observation.

% If the answer to question~\ref{Q2} is no, it means that the method we implemented for $N > 40$ can not
% be complete. If the answer is yes, then the above degree tables are likely complete for many $N$'s
% (but probably not for every $N \leq 80$).
% There are infinitely many modular units, and we could only test a finite (randomly chosen) subset.
% So the search that produced the above tables is probabilistic, and there is a non-trivial chance that
% some degree(s) is/are missing, a chance that increases with $N$ and the rank of the Jacobian.

\end{document}